\documentclass[twoside,psfig]{article}

\usepackage[dvips]{graphicx}
\usepackage{amsmath}
\usepackage{amsfonts}
\usepackage{amssymb}
\usepackage{amscd}
\usepackage{amsthm}
\usepackage[all]{xy}

\oddsidemargin=\evensidemargin 
\catcode`\@=11 \long\def\@makefntext#1{ \protect\noindent \hbox to
3.2pt {\hskip-.9pt
$^{{\eightrm\@thefnmark}}$\hfil}#1\hfill}       

\def\ps@myheadings{\let\@mkboth\@gobbletwo      
\def\@oddhead{\hbox{}
\rightmark\hfil\eightrm\thepage}
\def\@oddfoot{}\def\@evenhead{\eightrm\thepage\hfil
\leftmark\hbox{}}\def\@evenfoot{}
\def\sectionmark##1{}\def\subsectionmark##1{}}

\catcode`@=11                        
\def\ps@plain{\let\@mkboth\@gobbletwo
     \def\@oddhead{}\def\@oddfoot{\eightrm\hfil\thepage
     \hfil}\def\@evenhead{}\let\@evenfoot\@oddfoot}

\newcounter{sectionc}\newcounter{subsectionc}\newcounter{subsubsectionc}
\renewcommand{\section}[1] {\vspace{12pt}\addtocounter{sectionc}{1}
\setcounter{equation}{0}
\setcounter{theorem}{0}\setcounter{lemma}{0}\setcounter{proposition}{0}
\setcounter{example}{0}\setcounter{remark}{0}\setcounter{corollary}{0}
\setcounter{subsectionc}{0}\setcounter{subsubsectionc}{0}\noindent
    {\tenbf\thesectionc. #1}\par\vspace{5pt}}
\renewcommand{\subsection}[1] {\vspace{12pt}\addtocounter{subsectionc}{1}
    \setcounter{subsubsectionc}{0}\noindent
    {\bf\thesectionc.\thesubsectionc.
    {\kern1pt \bfit #1}}\par\vspace{5pt}}
\renewcommand{\subsubsection}[1] {\vspace{12pt}
    \addtocounter{subsubsectionc}{1}
    \noindent
    {\tenrm\thesectionc.\thesubsectionc.\thesubsubsectionc. {\kern1pt
    \it #1}}\par\vspace{5pt}}
\newcommand{\nonumsection}[1] {\vspace{12pt}\noindent{\tenbf #1}
    \par\vspace{5pt}}

\newtheorem{theorem}{Theorem}[sectionc]

\newtheorem{proposition}{Proposition}[sectionc]
\newtheorem{example}{Example}[sectionc]
\newtheorem{remark}{Remark}[sectionc]
\newtheorem{corollary}{Corollary}[sectionc]

\newenvironment{proofs}{\begin{trivlist}
    \item[\noindent]{\bf Proof.}}{\quad $\square$\end{trivlist}}

\newcommand{\textlineskip}{\baselineskip=13pt}
\newcommand{\smalllineskip}{\baselineskip=10pt}

\newcommand{\copyrightheading}[1]
    {\vspace*{-2.5cm}\smalllineskip{\flushleft
    {\footnotesize Journal of Knot Theory and Its Ramifications #1}\\
    {\footnotesize \copyright\kern2pt World Scientific
         Publishing Company}\\
         }}

\def\abstracts#1#2#3#4{{
    \centering{\begin{minipage}{4.5in}\footnotesize\baselineskip=10pt
    \centerline{ABSTRACT}
    \parindent=15pt #1\par
    \parindent=15pt #2\par
    \parindent=15pt #3\par
    \parindent=15pt #4\par
    \end{minipage}}\par}}

\def\keywords#1{{
    \centering{\begin{minipage}{4.5in}\footnotesize\baselineskip=10pt
    {\footnotesize\it Keywords}\/: #1
    \end{minipage}}\par}}

\def\subject#1{{
    \centering{\begin{minipage}{4.5in}\footnotesize\baselineskip=10pt
    {\footnotesize  Mathematics Subject Classification 2000}\/: #1
    \end{minipage}}\par}}

\renewenvironment{thebibliography}[1]
    {\frenchspacing
     \ninerm\baselineskip=11pt
     \begin{list}{[\arabic{enumi}]}
    {\usecounter{enumi}\setlength{\parsep}{0pt}
     \setlength{\leftmargin 13.7pt}{\rightmargin 0pt} 
     \setlength{\itemsep}{0pt} \settowidth
    {\labelwidth}{[#1]}\sloppy}}{\end{list}}

\newcounter{itemlistc}
\newcounter{romanlistc}
\newcounter{alphlistc}
\newcounter{arabiclistc}

\newenvironment{romanlist}
    {\setcounter{romanlistc}{0}
     \begin{list}{$($\roman{romanlistc}$)$}
    {\usecounter{romanlistc}
     \setlength{\parsep}{0pt}
     \setlength{\itemsep}{0pt}}}{\end{list}}

\newcommand{\fcaption}[1]{
        \refstepcounter{figure}
        \setbox\@tempboxa = \hbox{\footnotesize Fig.~\thefigure. #1}
        \ifdim \wd\@tempboxa > 5in
           {\begin{center}
        \parbox{5in}{\footnotesize\smalllineskip Fig.~\thefigure. #1}
            \end{center}}
        \else
             {\begin{center}
             {\footnotesize Fig.~\thefigure. #1}
              \end{center}}
        \fi}

\newcommand{\tcaption}[1]{
        \refstepcounter{table}
        \setbox\@tempboxa = \hbox{\footnotesize Table~\thetable. #1}
        \ifdim \wd\@tempboxa > 5in
           {\begin{center}
        \parbox{5in}{\footnotesize\smalllineskip Table~\thetable. #1}
            \end{center}}
        \else
             {\begin{center}
             {\footnotesize Table~\thetable. #1}
              \end{center}}
        \fi}

\def\pmb#1{\setbox0=\hbox{#1}
    \kern-.025em\copy0\kern-\wd0
    \kern.05em\copy0\kern-\wd0
    \kern-.025em\raise.0433em\box0}

\def\fnt#1#2{\footnotetext{\kern-.3em
    {$^{\mbox{\scriptsize #1}}$}{#2}}}

\def\fpage#1{\begingroup
\voffset=.3in
\thispagestyle{empty}\begin{table}[b]\centerline{\footnotesize #1}
    \end{table}\endgroup}

\font\tenrm=cmr10  \font\tenbf=cmbx10
\font\bfit=cmbxti10 at 10pt \font\ninerm=cmr9 
 \font\eightrm=cmr8

\def\o{\omega}
\def\a{\alpha}

\def\bb{\gamma}

\def\wa{x}
\def\ee{\varepsilon}

\def\o{\omega}
\def\dd{\delta}

\def\ps{\varphi}

\def\A{\textup{A}}

\def\SS{\mathbf{S}}
\def\L{\textup{L}}
\def\F{\mathbb{F}}

\def\H{\textup{H}}
\def\G{\Gamma}
\def\F{\textup{F}}

\def\TO{\textup{Tors}}

\def\L{\textup{L}}

\def\Z{\mathbb{Z}}
\def\Q{\mathbb{Q}}

\def\gcd{\textup{gcd}}

\def\thefigure{\thesectionc.\arabic{figure}}
\begin{document}
\setlength{\textheight}{7.7truein}

\thispagestyle{empty} \setcounter{page}{1}


\vspace*{0.88truein}

\fpage{1}

\centerline{\bf THE ALEXANDER POLYNOMIAL OF
$\mathbf{(1,1)}$-KNOTS}

\vspace*{0.37truein}

\centerline{\footnotesize A. CATTABRIGA}

 \baselineskip=12pt \centerline{\footnotesize\it
 Mathematics  Department, University of Bologna,}
 \baselineskip=10pt
\centerline{\footnotesize\it
 P.zza di Porta S.Donato, 5, 40126 Bologna, Italy}
  \baselineskip=10pt
\centerline{\footnotesize\it cattabri@dm.unibo.it}

\vspace*{0.225truein}


\vspace*{0.21truein}

\abstracts{In this paper we investigate the  Alexander polynomial
of $(1,1)$-knots, which are knots lying in a 3-manifold with genus
one at most, admitting a particular decomposition. More precisely,
we study the  connections between the Alexander polynomial and a
polynomial associated to a cyclic presentation of the fundamental
group of an $n$-fold strongly-cyclic  covering branched over the
knot $K$, which we call the $n$-cyclic polynomial of $K$. In this
way, we generalize to all $(1,1)$-knots, with the only exception
of those lying in $\SS^2\times\SS^1$, a result obtained by
J.~Minkus for 2-bridge knots and extended by the author and
M.~Mulazzani  to the case of $(1,1)$-knots in $\SS^3$. As
corollaries some properties of the Alexander polynomial of knots
in $\SS^3$ are  extended to the case of $(1,1)$-knots  in lens
spaces.}{}{}{}

\vspace*{10pt}
\vspace*{10pt}
 \subject{ Primary
57M12, 57M27; Secondary  57N10.}

\vspace*{1pt}\textlineskip

\section{Introduction} \label{intro}

\begin{figure}[htbp]
\vspace*{13pt}
\begin{center}
\includegraphics*[totalheight=3.2cm]{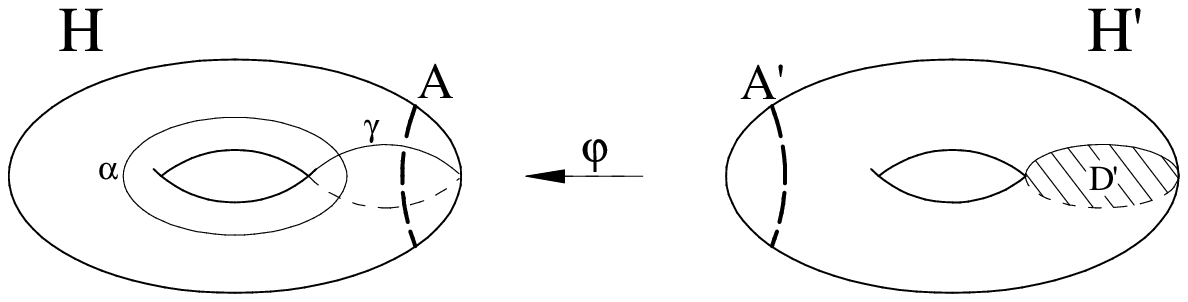}
\end{center}
\vspace*{13pt} \fcaption{A $(1,1)$-decomposition.} \label{figura}
\end{figure}

 \noindent A knot $K$ in a 3-manifold $M$ is called a
$(1,1)$-knot, or a 1-bridge torus knot, if it can be embedded into
a Heegaard torus of $M$ except at one over (or under) bridge. This
means that there exists a Heegaard splitting of genus one
\begin{equation}(M,K)=(\H,\A)\cup_{\ps}(\H',\A'),\end{equation}
where $\H$ and $\H'$ are solid tori, $\A\subset \H$ and
$\A'\subset \H'$ are properly embedded trivial arcs\footnote{This
means that there exists a disk $\textup{D}\subset \H$ (resp.
$\textup{D}'\subset \H'$) with $\A\cap \textup{D}=\A\cap\partial
\textup{D}=\A$ and $\partial \textup{D}-\A\subset\partial \H$
(resp. $\A'\cap \textup{D}'=\A'\cap\partial \textup{D}'=\A'$ and
$\partial \textup{D}'-\A'\subset\partial \H'$).}, and
$\ps:(\partial \H',\partial \A')\to(\partial \H,\partial \A)$ is
an attaching homeomorphism (see Fig. \ref{figura}). Such a
decomposition is called \textit{$(1,1)$-decomposition} of the
knot. By definition, the ambient manifold has genus less or equal
than one, so $M$, up to homeomorphism, is $\SS^2\times\SS^1$ or a
lens space $\L(p,q)$, including the case $\L(1,0)=\SS^3$.
Important examples of $(1,1)$-knots in $\SS^3$ are torus knots (by
definition) and 2-bridge knots (see \cite{KS}). This family of
knots has interesting features and has recently been
 studied by many authors from different points of view.
For instance, in \cite{CM2,CM4,CK}, different representations of
$(1,1)$-knots are given, while in \cite{AGM, CM1,CM3,GM,Mu} their
connections with manifolds with cyclically presented fundamental
groups are studied.

 In this paper we investigate the Alexander polynomial of
$(1,1)$-knots. This class of knots seems to have very interesting
features with regards to the Alexander polynomial. For example, in
\cite{Fu}, the author shows that it is possible to realize each
Alexander polynomial in $\SS^3$ using $(1,1)$-knots. Moreover, he
exhibits an infinite family of $(1,1)$-knots in $\SS^3$ with a
trivial Alexander polynomial. In this paper we  focus especially
on $(1,1)$-knots in lens spaces. From the results of \cite{Mu},
each $n$-fold strongly-cyclic branched covering of a (1,1)-knot
admits a cyclic presentation for the fundamental group and  a
polynomial associated to this  presentation. Moreover, if
$K\subset\L(p,q)$ and $n$ is coprime with $p$,  the covering is
unique, and so there is a natural way to associate to $K$ a
polynomial that will be called the $n$-cyclic polynomial of $K$.
Using the reduced Reidemeister torsion of the knot complement, we
find the relation between this polynomial and the Alexander
polynomial. In this way, we generalize to all $(1,1)$-knots, with
the only exception of those lying in $\SS^2\times\SS^1$, a result
of \cite{Mi} for 2-bridge knots and of \cite{CM1} for
$(1,1)$-knots in $\SS^3$. As corollaries, we extend to all
$(1,1)$-knots in lens spaces some properties of the Alexander
polynomial of  knots in $\SS^3$.\\

In Section 2 we  give the definition of $n$-cyclic polynomials. In
Section 3 we recall the definition of the Alexander polynomial of
a knot, and some results concerned with it.  The main result of
this paper is established in  Section 4, where  some corollaries
and examples are also collected.

\newpage
\section{Strongly-cyclic coverings and cyclic presentations of groups}
 \noindent Let $\F_n$ denote the free group of rank $n$. A
group $G$ is called  \textit{cyclically presented} if there exists
$n>0$ and $w\in\F_n$ such that \begin{equation}G=G_n(w)=\langle
x_1,\ldots,x_n\mid
w,\theta_n(w),\ldots,\theta_n^{n-1}(w)\rangle,\end{equation} where
$\theta_n:\F_n\to \F_n$ is the automorphism defined by $\theta_n
(x_i)=x_{i+1}$ (subscripts mod $n$), for $i=1,\ldots,n$. Such a
presentation is called a \textit{cyclic presentation}. Obviously
$G_n(w)$ has an automorphism of order $n$ induced by $\theta_n$.

The \textit{polynomial associated to the  cyclic presentation}
$G_n (w)$ is \begin{equation} p_{n,w}(t)= \sum_{i=1}^{n} a_i
t^{i-1},\end{equation} where $a_i$ is the total exponent sum of
the letter $x_i$ in the word $w$, for \hbox{$i=1,\ldots, n$}.

For each $s=0,\ldots,(n-1)$, we have $G_n(w)=G_n(\theta_n^s(w))$,
while the polynomials $p_{n,\theta_n^s(w)}(t)$  are different
elements of $\Z[t]$. However, their image  in the quotient ring
$\Z[t]/(t^n-1)$ is the same, up to units. So, from now on, we will
consider the polynomial  associated to the cyclic presentation
$G_n(w)$ as an element of $\Z[t]/(t^n-1)$.

Note  that the abelianization of $G_n(w)$ has a circulant
presentation matrix (as a $\Z$-module), whose first row is given
by
the coefficients of  $p_{n,w}$.\\

An $n$-fold cyclic covering $f:\widetilde{M}\to M$, branched over
a knot $K\subset M$, is called {\it strongly-cyclic} if the
branched index of $K$ is $n$. In other words\footnote{We denote
with $N(K)\subset M$ an open tubular neighborhood of $K$ in M.},
if \hbox{$\o_f:H_1(M-N(K))\to \Z_n$} denotes the monodromy of $f$
and $\overline{m}$ is the first homology class of a meridian loop
of $K$, $\o_f(\overline{m})$ generates $\Z_n$. Moreover, two
$n$-fold strongly-cyclic coverings $f$ and $f'$, branched over the
same knot $K\subset M$, are \hbox{equivalent} if there exists an
invertible element $u\in\Z_n$ such that the following diagram
commutes
\begin{equation}\xymatrix{ H_1(L(p,q)- N(K)) \ar[dd]_{\o_f}\ar[rr]^{\o_{f'}}
& & \Z_n \ar[ddll]^{\mu_u} \\ & & \\   \Z_n & &}
\end{equation}
where $\mu_u$ denotes the multiplication by $u$. So, up to
equivalence, a cyclic branched covering $f$ is strongly-cyclic if
$\o_f(\overline{m})=1$. From now on, $C_{n}(K)$ will denote an
$n$-fold strongly-cyclic branched covering of $K$.

If we are dealing with knots in $\SS^3$, every $n$-fold cyclic
branched covering of a knot  is strongly-cyclic, so it exists and
 is unique, up to equivalence. This is no longer true for a
$(1,1)$-knot $K$ in a lens space $\L(p,q)$. However, from the
existence and uniqueness conditions of \cite{CM1}, it follows that
for each $n$ coprime with $p$  there is a unique $n$-fold
strongly-cyclic branched covering of $K$.

The connection between $(1,1)$-knots and cyclically presented
groups is given by the following theorem.
\begin{theorem}\textup{\cite{Mu}}  Every $n$-fold strongly-cyclic branched covering of a $(1,1)$-knot
admits a cyclic presentation for the fundamental group with $n$
generators.
\end{theorem}
Moreover,  this presentation is obtained by lifting a
$(1,1)$-presentation of the knot, as follows (for details see
\cite{CM1}).

From a $(1,1)$-decomposition
$(\L(p,q),K)=(\H,\A)\cup_{\ps}(\H',\A')$ of  $K$,  we get the
following presentation of the knot group:
\begin{equation}
 \label{pi1} \pi_1(\L(p,q)-N(K),\ast)=\langle
\a,\bb,\,\mid\, r(\a,\bb)\rangle,
\end{equation} where the generators $\a$ and $\bb$ represent, respectively, a longitude
of $\partial\H$  and a meridian loop of $K$, and the relator
$r(\a,\bb)$ corresponds to the loop $\ps(\partial \textup{D}')$,
where $\textup{D}'$ is a meridian disk of $\H'$ that does not
intersect $\A'$ (see Fig. \ref{figura}).

If necessary, by replacing the covering
$f:\widetilde{M}\to\L(p,q)$ with an equivalent one, we can suppose
that $\o_f(\bb)=1$. Let $\bar r(\wa,\bb)=
r(\wa\bb^{\o_{f}(\a)},\bb)=\wa^{\ee_1}\bb^{\dd_1}\cdots\wa^{\ee_s}\bb^{\dd_s}$
for some $\ee_1,\ldots,\ee_s,\dd_1,\ldots,\dd_s\in\Z$. We have
\hbox{$\pi_1(\widetilde{M},\widetilde{*})=G_n(w)$} where:
\begin{equation}
\label{word}
 w=x_{i_1}^{\ee_1}\cdots
x_{i_s}^{\ee_s}\ \ \ \textup{ (subscripts mod n)}, \end{equation}
with $i_k\equiv 1+\sum_{j=1}^{k-1}\dd_j \mod n$, for
$k=1,\ldots,s$.

For each $n$ coprime with $p$, the polynomial associated to the
cyclic presentation of the (unique) $n$-fold strongly-cyclic
branched covering of $K$ obtained as above (i.e. lifting a
$(1,1)$-decomposition of $K$ and under  the condition
$\o_f(\bb)=1$) will be called the $n$\textit{-cyclic polynomial}
of $K$, and it will be denoted with $\G_{K,n}$.

\begin{remark}\textup {In the case of a $(1,1)$-knot $K\subset\SS^2\times\SS^1$, which corresponds to $p=0$, the
situation is rather different. We do not have the uniqueness of
the $n$-fold strongly-cyclic branched covering of $K$ for any
value of $n>1$, and so we do not have a natural way of defining
the $n$-cyclic polynomial. Moreover,  there exists at most a
finite number of  finite strongly-cyclic branched covering of $K$,
for almost all the $K$.}\end{remark}

Let $\Delta_K\in\Z[t,t^{-1}]$ be the Alexander polynomial of a
knot $K\subset\SS^3$ and denote with $\Delta_{K,n}$ its projection
on $\Z[t]/(t^n-1)$. With these notations the results of \cite{CM1}
and \cite{Mi} can be restated as follows.
\begin{proposition} \label{gen} Let $K\subset\SS^3$.  The following holds:
\begin{romanlist} \item[\textup{(i)}] \textup{\cite{Mi}} if $K$ is a
2-bridge knot, for each $n>1$, we have
$\G_{K,n}(t)=\Delta_{K,n}(t)$, up to units of $\Z[t]/(t^n-1)$.
\item[\textup{(ii)}] \textup{\cite{CM1}} if $K$ is a $(1,1)$-knot,
for each $n>1$, we have $\G_{K,n}(t)=\Delta_{K,n}(t)$, up to units
of $\Z[t]/(t^n-1)$.
\end{romanlist}
\end{proposition}

The main result of this article is the generalization of this
relation to all  $(1,1)$-knots in lens spaces.

\section{Alexander polynomial}
\noindent In this section we recall the definition of the
Alexander polynomial of a knot $K$ in a compact, connected
3-manifold, and give some of its characteristics.

Let $M$ be a compact connected manifold, $\ast\in M$ be a fixed
point and denote with $\Phi:\pi_1(M,\ast)\to H_1(M)$ the Hurewitz
homomorphism. Consider the projection $j:H_1(M)\to
H_1(M)/\TO(H_1(M))$ and the induced ring homomorphism\footnote{For
a group $G$  we denote with $\Z[G]$ its  integral group ring.}
$\bar{j} :\Z[H_1(M)]\to \Z[H_1(M)/\TO(H_1(M))]$. Note that, if $r$
is the first Betti number of $M$, we have
 $\Z[H_1(M)/\TO(H_1(M))]\cong\Z[t_1,t_1^{-1},\ldots ,t_r,t_r^{-1}]$, where    $t_1,\ldots,t_r$ are
generators of $H_1(M)/\TO(H_1(M))$. Let $E_1(M)\subset \Z[H_1(M)]$
be the first elementary ideal of $\pi_1(M,*)$ (see \cite{CF}) and
denote with $\bar{E}_1(M)$ the smallest principal ideal of
$\Z[H_1(M)/\TO(H_1(M))]$ containing $\bar{j}(E_1(M))$. The
generator $\Delta_M$ of $\bar{E}_1(M)$ is well-defined up to
multiplication by units of $\Z[H_1(M)/\TO(H_1(M))]$  and is called
the \textit{Alexander polynomial} of $M$.

If $K\subset X$ is a knot in a compact connected 3-manifold, the
Alexander polynomial of $K$ is the Alexander polynomial of
$M=X-N(K)$.

For a knot $K\subset\SS^3$ the Alexander polynomial has the
following characterization.

\begin{theorem} \label{gruppo}\textup{\cite{CF}} Let $\Delta_K(t)\in\Z[t,t^{-1}]$ be the Alexander
polynomial of a knot $K\subset\SS^3$. We have:
\begin{romanlist}\item[\textup{(i)}] $\Delta_K(t)=\Delta_K(t^{-1})$.
\item[\textup{(ii)}] $\Delta_K(1)=\pm 1$.
\end{romanlist}
Moreover, each polynomial $q(t)\in\Z[t,t^{-1}]$ satisfying
conditions 1 and 2 is the Alexander polynomial of a knot in
$\SS^3$.
\end{theorem}
Actually, from \cite{Fu}, each Alexander polynomial of a knot in
$\SS^3$ can be realized by a $(1,1)$-knot in $\SS^3$.

Moreover, the Alexander polynomial of a knot $K\subset\SS^3$
determines the order of the first homology group of $C_n(K)$, the
$n$-fold (strongly-)cyclic branched covering of $K$.

\begin{theorem}\textup{\cite{FO}} The abelian group $H_1(C_n(K))$ is finite if and
only if no root of the Alexander polynomial $\Delta_K(t)$ of $K$
is an $n$-root of unity. In this case: \begin{equation}\sharp
H_1(C_n(K))=\vert\prod_{\zeta^n=1}
\Delta_K(\zeta)\vert.\end{equation}
\end{theorem}

Let $K\subset\L(p,q)$ be a $(1,1)$-knot. Then
\hbox{$H_1(\L(p,q)-N(K))\cong\Z\oplus\Z_d$}, where $d>0$ divides
$p$ (see \cite{CM1}).  Then the  first Betti number of
$\L(p,q)-N(K)$ is one, so $\Delta_K$ lies in $\Z[t,t^{-1}]$.
Moreover, it is easy to check that \hbox{$\chi(\L(p,q)-N(K))=0$}.
So the following theorem holds with \hbox{$M=\L(p,q)-N(K)$}.

\begin{theorem}\label{Tu}\textup{\cite[Theorem
B]{T}} Let $M$ be  a compact, connected, orientable
\hbox{3-manifold} with $\chi(M)=0$. If $\partial M\ne\emptyset$
and $b_1(M)=1$, then  $r_M(t)=\Delta_M(t)/(t-1)$, where
$r_M(t)\in\Q(t)$ denotes the reduced Reidemeister torsion of $M$
and $t$ is a \hbox{generator} of $H_1(M)/\TO(H_1(M))$.
\end{theorem}

From  this theorem it follows (see \cite{T,T1}) that   for the
Alexander polynomial of $(1,1)$-knots in lens spaces it is still
true that $\Delta_K(t)=\Delta_K(t^{-1})$. However, we will see
that the property $\Delta_K(1)=\pm 1$ is no longer true.

For details on the reduced Reidemeister torsion, and on its
calculation,  we refer to \cite{T,T1}.

\begin{remark} \textup{ In the case of a $(1,1)$-knot
$K\subset\SS^2\times\SS^1$ it could happen that the first Betti
number is two. That is the case, for example, of the trivial
knot.}
\end{remark}

\section{Main theorem}

\begin{theorem} \label{Main}Let $K\subset\L(p,q)$ be a
$(1,1)$-knot and denote with $d$ the order of the torsion subgroup
of $H_1(\L(p,q)-N(K))$. Then, for each $n>1$ such that
\hbox{$\gcd(p,n)=1$,} we have:
\begin{equation}\G_{K,n}(t^{\frac{p}{d}})=\Delta_{K,n}(t)\sum_{i=0}^{\frac{p}{d}-1}t^i,\end{equation}
up to units of $\Z[t]/(t^n-1)$, where $\G_{K,n}$ is the $n$-cyclic
polynomial associated to $K$ and $\Delta_{K,n}$ is the projection
of the Alexander polynomial of $K$ on $\Z[t]/(t^n-1)$.
\end{theorem}
\begin{proofs}   To prove the statement we will use the reduced Reidemeister torsion of
\hbox{$M=\L(p,q)-N(K)$.} In the cellular presentation (\ref{pi1})
of $\pi_1(M,*)$, the total exponent sum of $\a$ in $r(\a,\bb)$ is
$p$, since $\a$
 is a generator of $\pi_1(\L(p,q),\ast)$. So, by abelianization,  we get $H_1(M)=\langle
\a,\bb\,\mid\,p\a+q'\bb\rangle\cong\Z\oplus\Z_d$, where
$d=\gcd(p,q')$. Moreover, if we set $\bar{p}=p/d$ and
$\bar{q}=q'/d$ we have that $\eta=\bar{p}\a+\bar{q}\bb$ and
$\xi=-s\a+r\bb$ generate, respectively, the torsion part and the
free part of $H_1(M)$, where $r$ and $s$ are integers such that
$r\bar{p}+s\bar{q}=1$. Now,  consider the covering
$F:\overline{M}\to M$ corresponding to the subgroup
$\Phi^{-1}(\TO(H_1(M)))$ of $\pi_1(M,*)$, where
$\Phi:\pi_1(M,\ast)\to H_1(M)$ denotes the Hurewitz homomorphism.
The  Betti number of $H_1(M)$ is one, so the covering $F$ is
infinite cyclic with monodromy $\omega_F:H_1(M)\to\langle
t\rangle\cong\Z$ defined by $\omega_F(\xi)=t$ and
$\omega_F(\eta)=1$. Inverting the defining relations for $\xi$ and
$\eta$ we have that $\omega_F(\bb)=t^{\bar{p}}$ and
$\omega_F(\a)=t^{-\bar{q}}$. Up to contraction, which does not
alter the torsion,  the cell complex $C_*(M)$
\begin{equation}0\to C_2=\langle \ps(\textup{D}')\rangle\to
C_1=\langle\a,\bb\rangle\to C_0=\langle \ast\rangle\to
0,\end{equation} is a cellular decomposition for $M$. The reduced
Reidemeister torsion of $M$ is, up to multiplication by $\pm
t^{\pm h}$, the torsion of the complex
$\overline{C}_*(\overline{M})=\Q(t)\otimes_{\Z[t,t^{-1}]}C_*(\overline{M})$,
where $C_*(\overline{M})$ is the lifting of $C_*(M)$, with  a
fundamental family of cells of $\overline{M}$ as a base. Let
$\delta_1:\overline{C_1}\to\overline{C_0}$ and
$\delta_2:\overline{C_2}\to\overline{C_1}$ be the boundary
operators. If $\tilde{*}$ is the fixed 0-cell over $*$, and
$\tilde{\a}$ and $\tilde{\bb}$ are the lifting of, respectively,
$\a$ and $\bb$ with starting point $\tilde{*}$ then, by the action
of the monodromy map, we have
$\delta_1(\tilde{\a})=(t^{-\bar{q}}-1)\tilde{*}$ and
$\delta_1(\tilde{\bb})=(t^{\bar{p}}-1)\tilde{*}$. Moreover it is
well-known (see for example \cite[$\S$ 9A-9B]{BZ}) that, the
matrix of $\delta_2$, is the Alexander matrix of the  presentation
(\ref{pi1}). So $\delta_1$ and $\delta_2$ are represented by the
matrices $^\textup{t}(t^{-\bar{q}}-1\ \ \ t^{\bar{p}}-1)$ and
$(Q_{\a}(t)\ \ \ Q_{\bb}(t))$, where $Q_{\a}(t)$ (resp.
$Q_{\bb}(t)$) is obtained from the Fox derivative $\partial
r(\a,\bb)/\partial\a$ (resp. $\partial r(\a,\bb)/\partial\bb$) by
the substitutions $\a=t^{-\bar{q}}$ and $\bb=t^{\bar{p}}$. The
collection $\{1,1,\emptyset; t^{\bar{p}}-1,Q_{\a}(t)\}$ is a
torsion chain for this complex (see \cite{T,T1}), so
$r_M(t)=Q_{\a}(t)/(t^{\bar{p}}-1)\in\Q(t)$. By Theorem \ref{Tu},
$r_M(t)=\Delta_K(t)/(t-1)$ and so
$Q_{\a}(t)=\Delta_K(t)(\sum_{i=0}^{\bar{p}-1}t^i)$, where the
equality holds up to units of $\Z[t,t^{-1}]$. To complete the
proof we have to show that $Q_{\a}(t)=\G_{K,n}(t^{\bar{p}})$ in
$\Z[t]/(t^n-1)$. Observe that, for each $n>1$ coprime with $p$,
the monodromy of the unique $n$-fold strongly-cyclic branched
covering of $K$, $F_n$, is the composition of $\omega_F$ with the
epimorphism $\Z\to\Z_n$, given by $t\to 1$. Then we have
$\o_{F_n}(\a)=-\bar{q}$ and $\o_{F_n}(\bb)=\bar{p}$. If
$\bar{p}\ne 1$, to  calculate $\G_{K,n}$ we have to replace $F_n$
with the equivalent covering with monodromy $\o'$, such that
$\o'(\bb)=1$, and so $\o'(\a)=-\bar{q}\bar{p}^{-1}$, where
$\bar{p}^{-1}$ is the inverse of $\bar{p}$ in $\Z_n$.  Let
$\bar{r}(x,\gamma)=r(x\gamma^{\omega'(\a)},\bb)$. It is easy to
check, using formula (\ref{word}), that $\G_{K,n}(u)$ is equal to
the polynomial obtained
 from $\partial \bar{r}(x,\bb)/\partial x$ by the
substitutions $x=1$ and $\bb=u$. Moreover, since $\partial
\bar{r}(x,\bb)/\partial x=\partial r(x\bb^{\o'(\a)},\bb)/\partial
x=(\partial r/\partial \a)(x\bb^{\o'(\a)},\bb)$ we have that
$\G_{K,n}(u)$ is equal to the polynomial obtained from $\partial
r/\partial \a$ by the substitutions \hbox{$
\a=u^{\o'(\a)}=u^{-\bar{q}\bar{p}^{-1}}$,} and $\bb=u$. So,
setting $u=t^{\bar{p}}$ in $\G_{K,n}(u)$, we get $Q_{\a}(t)$,
which ends the proof.
\end{proofs}
\begin{remark} \textup{In the case of
$\SS^2\times\SS^1$, the previous theorem does not hold, since, as
we have already observed, in this case $p=0$, so for each $n>1$,
$\gcd(n,p)=n\ne 1$. In fact, if we look at the proof, it is easy
to see that, since \hbox{$\bar{p}=p/d=0$},  the monodromy
$\omega_F$ of the covering  that leads to the calculation of the
Alexander polynomial sends a meridian of the knot into $t^0$. So,
the covering whose monodromy is the composition of $\omega_F$ with
the projection $\Z\to\Z_n$, defined by $t\to 1$, is not
strongly-cyclic for any value of $n>1$.}
\end{remark}
Observe that, when  $p=1$, $K$ is  a $(1,1)$-knot in $\SS^3$, and
we get the same statement as in Proposition \ref{gen}. Moreover,
if $\TO(H_1(\L(p,q)-N(K)))\cong\Z_p$, then the  projection of the
Alexander polynomial on $\Z[t]/(t^n-1)$ is equal to the $n$-cyclic
polynomial, up to units.

As a corollary of this theorem, we get a generalization of Theorem
\ref{gruppo}.

\begin{corollary}Let $K\subset\L(p,q)$ be a $(1,1)$-knot, and let $C_n(K)$ be the \hbox{$n$-fold} strongly-cyclic branched covering of $K$,
with  $n$ coprime with  $p$. Denote with $d$ the order of the
torsion subgroup of $H_1(\L(p,q)-N(K))$. Then $H_1(C_n(K))$ is
finite if and only if no root of $\Delta_{K,n}$  is an $n$-root of
unity. Moreover if  $d_n$ denote the order of the  torsion
subgroup of $H_1(C_n(K))$, we have:
\begin{equation}
d_n=\vert\prod_{\zeta^n=1,\Delta_{K,n}(\zeta)\ne
0}\left(\frac{\Delta_{K,n}}{\Phi}\right)(\zeta)\sum_{j=0}^{\frac{p}{d}-1}\zeta^j\vert,\end{equation}
where $\Phi$ is the product of the distinct cyclotomic polynomials
$\Phi_s$ such that $s$ divides $n$ and $\Phi_s$ divides
$\Delta_{K,n}$. In particular, if $H_1(C_n(K))$ is finite, we
have:
\begin{equation}\sharp
H_1(C_n(K))=\vert\prod_{\zeta^n=1}\Delta_{K,n}(\zeta)\sum_{j=0}^{\frac{p}{d}-1}\zeta^j\vert.\end{equation}

\end{corollary}

\begin{proofs}
If  $\o$ denotes the monodromy of  $C_n(K)$, up to equivalence, we
can suppose that  $\o(\bb)=1$. Consider the equivalent $n$-fold
strongly-cyclic branched covering $C'_n(K)$ with monodromy
$\o'=\mu_{p/d}\o$, where $\mu_{p/d}:\Z_n\to \Z_n$ denotes the
multiplication by $p/d$ which is invertible in $\Z_n$. Then the
isomorphism \hbox{$H_1(C_n(K))\to H_1(C'_n(K))$} is given by
$x_i\to x_j$ where $j\equiv i\dot(p/d)\ \mod\ n$. As previously
observed, the circulant matrix whose first row is given by the
coefficients of $\G_{K,n}(t)$ is a  presentation matrix for
$H_1(C_n(K))$ as a $\Z$-module, so, the circulant matrix $B$ whose
first row is given by the coefficients of $G(t)=\G_{K,n}(t^{p/d})$
is a presentation matrix for $H_1(C'_n(K))$ as a $\Z$-module.
Obviously $H_1(C'_n(K))$ is finite if and only if $H_1(C_n(K))$ is
finite and the order of the torsion subgroup of $H_1(C'_n(K))$ is
$d_n$. By the theory of circulant matrices (see \cite{BM}), there
exists a complex unitary matrix $F$, called Fourier matrix, such
that
$FBF^*=D=\textup{Diag}(G(\zeta_1),G(\zeta_2),\ldots,G(\zeta_n))$,
where $\zeta_1,\zeta_2,\ldots,\zeta_n$ are the $n$-roots of the
unity. So $H_1(C'_n(K))$ is finite if and only if the rank of $B$
is $n$, and so if and only if
\hbox{$G(\zeta_i)=\Delta_{K,n}(\zeta_i)\sum_{j=0}^{\frac{p}{d}-1}\zeta_i^j\ne
0$} for each \hbox{$i=1,\ldots,n$}. The first statement follows
from the fact that, since $\gcd(p,n)=1$, we have $\gcd(p/d,n)=1$,
and so $q_{p/d}(\zeta_i)\ne 0$, for $i=1,\ldots,n$, where
$q_{p/d}(t)=\sum_{j=0}^{\frac{p}{d}-1}t^j$. Moreover, by
\cite[Theorem 3.3]{SW}, we have that
$d_n=\vert\prod_{\zeta^n=1,G(\zeta)\ne
0}\left(\frac{G}{\Phi}\right)(\zeta)\vert,$ where $\Phi$ is the
product of the distinct cyclotomic polynomials $\Phi_s$ such that
$s$ divides $n$ and $\Phi_s$ divides $G$.  To end the proof, we
have only to show that $\Phi_s$ is a divisor of
$\Delta_{K,n}q_{p/d}$ if and only if  is a divisor of
$\Delta_{K,n}$. This follows from the fact that, if  $s$ is a
divisor  of $n$, we have $\gcd(s,p/d)$=1, so for each $\zeta$
primitive  $s$-root of the unity, $q_{p/d}(\zeta)\ne 0$.
\end{proofs}

Another straightforward corollary is the following.

\begin{corollary}The Alexander polynomial of a $(1,1)$-knot $K\subset\L(p,q)$ satisfies the relation $\Delta_K(1)=\pm
d$, where $d$ is the order of the torsion subgroup of
\hbox{$H_1(\L(p,q)-N(K))$}.
\end{corollary}
\begin{proofs} It is enough to observe that $\G_{K,n}(1)=\pm p$.
\end{proofs}
An interesting question for future study could be whether any
symmetric polynomial in $\Z[t,t^{-1}]$ can be realized as the
Alexander polynomial of a $(1,1)$-knot in a lens space.
\begin{example} \textup{Let  $K_{p,q}$ be the trivial knot in $\L(p,q)$.
We have that \hbox{$\pi_1(\L(p,q)-N(K_{p,q}),\ast)=\langle
\a,\bb\,\mid\,\a^p\rangle$,} (for reference see \cite{CM2}). Then
\hbox{$p=d=\sharp\TO(H_1(\L(p,q)-N(K_{p,q})))$.} So, for each
$n>p$ with $\gcd(n,p)=1$, we have
$\G_{K_{p,q},n}=p=\Delta_{K_{p,q}}.$}
\end{example}

We end this paper by observing that formula (\ref{word}) can
easily be implemented to find $\G_{K,n}$, and so, by Theorem
\ref{Main}, to calculate $\Delta_K$, as shown in the following
example.
\begin{example}
\textup{By results of \cite{CM4}, each $(1,1)$-knot $K$ can be
represented as $K(a,b,c,r)$, where $a,b,c,r$ are nonnegative
integer parameters  determining the Heegaard diagram of a
$(1,1)$-decomposition of $K$. For each $m>2$, let \hbox{$K_m
=K(1,m-2,0,1)\subset\L(m-2,1)$.} This family of knots is very
interesting since, as proved in  \cite{GM2}, the $m$-fold
strongly-cyclic branched covering of $K_m$, with \hbox{monodromy}
that sends $\bb$ into $1$ and $\a$ into $0$, is the Neuwirth
manifold $\mathcal{N}_m$ of type $m$. These manifolds were
introduced by L. Neuwirth in \cite{N}, while in \cite{C} it is
proved that $\mathcal{N}_m$ is  a Seifert manifold of type
$(0;-1;(2,1),\ldots, (2,1))$, with base $\SS^2$, Euler number $-1$
and   $m$ \hbox{exceptional} fibers. From the parametric
representation of $K_m$ (see \cite{AGM} for details) we get
\begin{equation}\label{ex}\pi_1(\L(m-2,1)-N(K_m),\ast)=\langle
\a,\bb\,\mid\,(\a\bb)^{m-1}\a^{-1}\bb\rangle, \end{equation} and
so
$$H_1(\L(m-2,1)-N(K_m))= \langle
\a,\bb\,\mid\,(m-2)\a+m\bb\rangle\cong\left\{
\begin{array}{ll}
    \Z & \hbox{if $m$ is odd,} \\
    \Z\oplus\Z_2 & \hbox{is $m$ is even.} \\
\end{array}
\right. $$ Notice that  $K_3$ is the  trefoil knot in $\SS^3$.}

\textup{If we denote with $\o$ the monodromy of the $n$-fold
strongly-cyclic branched covering of $K_m$, with $\gcd(n,m-2)$=1
and satisfying the condition $\o(\bb)=1$, from the homology
relation we get
$$c=\o(\a)=\left\{
\begin{array}{ll}
    -m\cdot\overline{m-2}, & \hbox{if $m$ is odd,} \\
    -\frac{m}{2}\cdot\overline{(\frac{m-2}{2})}, & \hbox{if $m$ is even,} \\
\end{array}\right.  $$
where $\overline{x}$ denote the inverse of $x$ in $\Z_n$. From the
substitution $x=\a\bb^{-c}$ in the presentation (\ref{ex}),  we
get the relator
\hbox{$\bar{r}(x,\bb)=(x\bb^{c+1})^{m-1}\bb^{-c}x^{-1}\bb$.} By
formula (\ref{word}), $\pi_1(C_n(K_m),\tilde{*})=G_n(w)$, where
$$w=(\prod_{i=0}^{m-2}x_{1+i(1+c)})x^{-1}_{2+(m-2)(1+c)}$$
(subscripts  $\mod\ n$). So, we get:
$$\G_{K_m,n}(t)=(\sum_{i=0}^{m-2}t^{i(1+c)})-t^{1+(m-2)(1+c)}.$$
If $m$ is odd then $p/d=m-2$ and we get
 $$\G_{K_m,n}(t^{m-2})=(\sum_{i=0}^{m-2}t^{-2i})-t^{-(m-2)}=(\sum_{i=0}^{m-2}t^{2i})-t^{(m-2)}=\Delta_{K_m,n}(t)\sum_{j=0}^{m-3}t^j,$$
(where the second equality holds up to multiplication for $t^{2(m-2)}$). This means that for each $n>2(m-2)$ coprime with $m-2$,
 the representative of $\Delta_{K_m}$ in $\Z[t]/(t^n-1)$ does not
 depend on $n$. So, for each $m>2$ odd
 $$\Delta_{K_m}(t)=\frac{(\sum_{i=0}^{m-2}t^{2i})-t^{(m-2)}}{\sum_{j=0}^{m-3}t^j}=\sum_{i=0}^{m-1}
 (-t)^i.$$
 Analogously, if $m$ is even we have $p/d=\frac{m-2}{2}$ and
 therefore
 $$\G_{K_m,n}(t^{(m-2)/2})=(\sum_{i=0}^{m-2}t^{-i})-t^{-(m-2)/2}=(\sum_{i=0}^{m-2}t^{i})-t^{(m-2)/2}=\Delta_{K_m,n}(t)\sum_{j=0}^{(m-4)/2}t^j.$$
By the same considerations, for each $m>2$ even
$$\Delta_{K_m}(t)=\frac{(\sum_{i=0}^{m-2}t^{i})-t^{(m-2)/2}}{\sum_{j=0}^{(m-4)/2}t^j}=t^{m/2}+1.$$}
\end{example}

\nonumsection{Acknowledgements}

 \noindent Work performed under the auspices of the G.N.S.A.G.A.
of I.N.d.A.M. (Italy) and the University of Bologna, funds for
selected research topics. The author  would like to thank Michele
Mulazzani and the referee for their helpful suggestions.

\nonumsection{References}

\end{document}